\documentclass[12pt,oneside,reqno]{article}
\usepackage{braket, amsxtra, xr}
\usepackage[utf8]{inputenc}

\usepackage{slashed}
\usepackage{graphicx}
\usepackage{xcolor}
\usepackage{pdfrender}
\pdfrender{TextRenderingMode=2,LineWidth=.5pt}
\usepackage{amsmath}
\usepackage{amsthm}
\usepackage{indentfirst}
\usepackage{amsfonts}
\usepackage{marginnote}
\usepackage{xcolor}
\usepackage{hyperref}
\usepackage{OldStandard}
\usepackage{multicol}
\usepackage{lipsum}
\usepackage{newtxmath}

\usepackage{geometry}

\usepackage{blindtext}
\usepackage{authblk}
\usepackage{mathtools}

\usepackage{xcolor}
\usepackage{pdfrender}
\pdfrender{TextRenderingMode=2,LineWidth=.5pt}

\usepackage{lipsum} 

\newtheoremstyle{plainindent}
  {\topsep}   
  {\topsep}   
  {\itshape}  
  {\parindent}
  {\bfseries} 
  {.}         
  {5pt plus 1pt minus 1pt} 
  {}          

\theoremstyle{plainindent}

\parskip=\baselineskip

\newcommand\myshade{85}
\colorlet{mylinkcolor}{blue}
\colorlet{mycitecolor}{red}
\hypersetup{linkcolor  = mylinkcolor,
	citecolor  = mycitecolor,
	urlcolor   = myurlcolor!\myshade!black,
	colorlinks = true,}
	
\title{\textsc{ on finite mellin transform via ramanujan's master theorem}}
\author{\textit{By}\textsc{ Omprakash Atale${}^{\dagger}$ }\footnote{${}^{\dagger}$E-mail: atale.om@outlook.com\\ Keywords: Mellin transform, Ramanujan's Master Theorem, Incomplete gamma function}}
\affil{\small ${}^{\dagger}$\textit{Department of Mathematics, Savitribai Phule Pune University,}\\ \textit{Pune-411001, India}}
\date{[10 September 2024]}
\begin{document}
\maketitle
\begin{abstract}
   This paper aims to show that by making use of Ramanujan's Master Theorem and the properties of the lower incomplete gamma function, it is possible to construct a finite Mellin transform for the function $f(x)$ that has infinite series expansions in positive integral powers of $x$. Some applications are discussed by evaluating certain definite integrals. The obtained solutions are also compared with results from Mathematica to test the validity of the calculations.
\end{abstract}
\newpage
\fontsize{11.5}{16}\selectfont
\begin{center}
    \textsc{\textbf{\large \S I. Introduction}}
\end{center}
Ramanujan in his quarterly reports \cite{1} derived the following theorem. If $f$ has expansion of the form
\begin{equation}
f(x)=\sum_{n=0}^{\infty}(-1)^{n} \frac{\phi(n)}{n !} x^{n}\tag{1.1} 
\end{equation}
where $\phi(n)$ has a natural and continuous extension such that $\phi(0) \neq 0$, then for $n>0$, we have
\begin{equation}
\int_{0}^{\infty} x^{n-1}\left(\sum_{n=0}^{\infty}(-1)^{n} \frac{\phi(n)}{n !} x^{n}\right) dx =\phi(-n)\Gamma(n),\tag{1.2}   
\end{equation}
where $n$ is any positive integer. Ramanujan's method for deriving his master's theorem was unconventional and his theorem had a problem with convergence of the integral. 

For Eqn. (1.2) to work, $f(x)$ should be finite and continuous on the interval $[0,\infty]$ but not necessarily at $0$ and $\infty$ and $x^n f(x)$ should vanish when $x$ becomes infinite.

Hardy established some boundaries to the value of $\phi$ and derived a theorem that is in all respects convergent. Following is Hardy's version of the above theorem \cite{2}. Let $\varphi(z)$ be an analytic (single-valued) function, defined on a half-plane $H(\delta)=\{z \in {C}: \Re (z) \geq-\delta\}$ for some $0<\delta<1 .$ Suppose that, for some $A<\pi, \phi$ satisfies the growth condition $|\phi(v+i w)|<C e^{P v+A|w|} \forall z=v+i w \in H(\delta)$. Let $0<x<e^{-p}$ the growth condition shows that the series $\Phi(x)=\phi(0)-x \phi(1)+x^{2} \phi(3) \ldots$ converges. The residue theorem yields
\begin{equation}
\Phi(x)=\frac{1}{2 \pi i} \int_{c-i \infty}^{c+i \infty} \frac{\pi}{\sin s \pi} \phi(-s) x^{-s} d s\tag{1.3}
\end{equation}
for any $0<c<\delta$. Observe that $\pi / \sin \pi s$ has poles at $s=-n$ for $n=0,1,2 \ldots$ with residue $(-1)^{n}$. The above integral converges absolutely and uniformly for $c \in(a, b)$ and $0<a<b<\delta$. Using Mellin inversion formula, $\forall 0<\Re s<\delta$, we get
\begin{equation}
\int_{0}^{\infty} x^{s-1}\left\{\phi(0)-x \phi(1)+x^{2} \phi(3) \ldots\right\} d x=\frac{\pi}{\sin s \pi} \phi(-s). \tag{1.4}
\end{equation}
The substitution $\phi(u) \rightarrow \phi(u) / \mathrm{\Gamma}(u+1)$ in Eqn. (1.4) establishes Ramanujan's master theorem in its original form (Eqn. (1.2)). For more on Hardy's class and its extension, the reader can refer to \cite{5}. And for more analogues and applications of Ramanujan's Master Theorem, the reader can refer to \cite{2,3,4,6,7, 14}.

\begin{center}
    {\large\textsc{\textbf{\S II. Finite Mellin Transform }}}
\end{center}

Let $s\in\mathbb{C}$ such that $\Re(s)>0$ and $t\in\mathbb{R}^{+}$. Then we define the lower incomplete gamma function as :
\begin{equation*}
    \gamma(s,t):=\int_{0}^{t}x^{s-1}e^{-x}dx.\tag{2.1}
\end{equation*}
Making use of the above integral representation, we state the following result.

\textbf{\textsc{Theorem 1}} \textit{ Let $\Re(s)>0$ and $t\in\mathbb{R}^{+}$, then}
  \begin{equation*}
      \int_{0}^{t}x^{s-1}\left(\sum_{n=0}^{\infty}\phi(n)\frac{(-x)^{n}}{n!}\right)dx=\gamma(s,t)\phi(-s).\tag{2.2}
  \end{equation*}
\textit{Proof:}  Replace $x$ with $mx$ in Eqn (2.2) to get
\begin{equation*}
    \gamma(s,t)m^{-s}:=\int_{0}^{t}x^{s-1}e^{-mx}dx\tag{2.3}
\end{equation*}
where $s,m>0$. Now let $m=r^k$ where $r>0$, multiply both sides by $f^{(k)}h^k/k!$ ($f$ shall be specified later) and sum on $k$, $0\leq k\leq \infty$, to obtain
\begin{equation*}
    \sum_{k=0}^{\infty}\frac{f^{(k)}h^k}{k!}\int_{0}^{t}e^{-r^kx}x^{s-1}dx=\gamma(s,t)\sum_{k=0}^{\infty}\frac{f^{(k)}(hr^{-n})^k}{k!}.\tag{2.4}
\end{equation*}
Expand the exponential in Maclaurin series, invert the order of summation and integration, and apply Taylor's theorem to deduce
\begin{equation*}
    \int_{0}^{t}x^{s-1}\sum_{n=0}^{\infty}\frac{f(hr^n+a)(-x^n)}{n!}dx=\gamma(s,t)f(hr^{-s}+a).\tag{2.5}
\end{equation*}
Now let $f(hr^{-s}+a)=\phi(-s)$, then Eqn.(2.5) can be rewritten as
\begin{equation*}
    \int_{0}^{t}x^{s-1}\sum_{n=0}^{\infty}\frac{\phi(n)(-x^n)}{n!}dx=\gamma(s,t)\phi(-s).\tag{2.6}
\end{equation*}
This completes our proof.

\textbf{Example 1.1.} Consider the following geometric series for $|x^k|<1$:
\begin{equation*}
    \frac{1}{1+x^k}=1-x^k+x^{2k}-x^{3k}+...=\sum_{n=0}^{\infty}(-x^k)^n.\tag{2.7}
\end{equation*}
Replacing $x$ with $x^k$, $s$ with $s/k (s\neq k)$ and substituting $\Phi(n)=\Gamma(n+1)$ in Theorem 1 gives
\begin{equation*}
      \int_{0}^{t}\frac{x^{s-1}}{1+x^k}dx=\gamma\left(\frac{s}{k},t\right)\Gamma\left(1-\frac{s}{k}\right).\tag{2.8}
  \end{equation*}

We may face limitations in simplifying the right-hand side of Eqn. (2.2) to simple expressions, unlike when using the gamma function instead of the lower incomplete gamma function, as in Ramanujan's Master Theorem (1.2). This is because the gamma function possesses more symmetries and properties than the lower incomplete gamma function. For instance, in Eqn. (1.2), when $\phi(n)=\Gamma(n+1)$, we can directly use the reflection formula for the gamma function to obtain an expression involving the sine function on the right-hand side of Eqn. (1.2). However, this is not possible in Example 1 because the lower incomplete gamma function does not have a corresponding reflection formula.

To test the validity of Eqn. (2.2), we evaluate for a particular case of $t=100, s=1/2, k=5$. In Mathematica, we get the following expression:
\begin{align*}
          &\int_{0}^{100}\frac{1}{\sqrt{x}(1+x^5)}dx=\bigg\{\frac{1}{5} \bigg(\sqrt[10]{-1} \bigg((-1)^{4/5} \log \left(1-10 \sqrt[10]{-1}\right)-(-1)^{4/5} \log \left(1+10
   \sqrt[10]{-1}\right)\\&+(-1)^{3/5} \log \left(1-10 (-1)^{3/10}\right)-(-1)^{3/5} \log \left(1+10
   (-1)^{3/10}\right)+\sqrt[5]{-1} \log \left(1-10 (-1)^{7/10}\right)-\\&\sqrt[5]{-1} \log \left(1+10
   (-1)^{7/10}\right)+\log \left(1-10 (-1)^{9/10}\right)-\log \left(1+10 (-1)^{9/10}\right)\bigg)+2 \tan
   ^{-1}(10)\bigg)\bigg\}.\tag{2.9}
\end{align*}
 Using  Eqn. (2.2) we get
\begin{equation*}
    \gamma\left(\frac{1}{10},100\right)\Gamma\left(\frac{9}{10}\right).
\tag{2.10}
\end{equation*}
It can be checked that both values are equivalent to $2.03328$.
\, 

\textbf{Example 1.2.} Consider the binomial series [\cite{1}, pg. 300]:
\begin{equation*}
    \frac{1}{(1+x)^k}=\sum_{n=0}^{\infty}\frac{\Gamma(n+k)}{\Gamma(k)}\frac{(-y)^n}{k!},\quad\,\,|y|<1,\tag{2.11}
\end{equation*}
we get $\phi(s)=\Gamma(s+k)/\Gamma(k)$. Using Theorem 1 we get the expression:
 \begin{equation*}
      \int_{0}^{t}\frac{x^{s-1}}{(1+x)^k}dx=\frac{\gamma(s,t)\Gamma(k-s)}{\Gamma(k)}.\tag{2.12}
  \end{equation*}
 
 We expect the above result to be restricted by some local constraints, such as $k\neq s$ and $\Re(s)>0$ since that would lead to an undefined value of the gamma function on the right-hand side of the above equation. These constraints would be different for different examples and would depend on the functions that appear on the right-hand side of the evaluated integral.
\, 

\textbf{Example 1.3.} Consider the following expansion,
\begin{equation*}
    \left(\frac{2}{1+\sqrt{1+4x}}\right)^\mu=\mu\sum_{k=1}^{\infty}\frac{\Gamma(2k+\mu)}{\Gamma(k+\mu+1)}\frac{(-x)^k}{k!}.\tag{2.13}
\end{equation*}
where $\mu$ is some natural number. By choosing $\phi(s)=\Gamma(2s+\mu)/\Gamma(s+\mu+1)$ in Theorem 1, we get
\begin{equation*}
      \int_{0}^{t}x^{s-1} \left(\frac{2}{1+\sqrt{1+4x}}\right)^\mu dx=\gamma(s,t)\frac{\Gamma(\mu-2s)}{(\mu-s+1)}\tag{2.14}
  \end{equation*}
  where $\mu-2s$ and $\mu-s+1$ should not be a negative integer. 

The following corollary follows directly from Theorem 1.
\, 

\textbf{\textsc{Corollary 1.1.}} \textit{ Let $\Re(s)>0$ and $t,u\in\mathbb{R}^{+}$, then}
  \begin{equation*}
      \int_{u}^{t}x^{s-1}\left(\sum_{n=0}^{\infty}\phi(n)\frac{(-x)^{n}}{n!}\right)dx=\left(\gamma(s,t)-\gamma(s,u)\right)\phi(-s).\tag{2.15}
  \end{equation*}
Unfortunately, such a simple formula is not known.

Now we make a small remark on how a Dirichlet series can be associated with Theorem 1. Replace $x$ with $mx$ in Eqn. (2.2), and multiply both sides with a Dirichlet coefficient $a_m$ to get the following result:
\begin{equation*}
      \int_{0}^{t}x^{s-1}\sum_{m=1}^{\infty}a_m f(mx)dx=\gamma(s,t)\phi(-s)g(s)\quad\mathrm{where}\quad g(s)=\sum_{m=1}^{\infty}\frac{a_m}{m^s}.\tag{2.16}
  \end{equation*}

There is no restriction to extend Theorem 1 to the definition of upper incomplete gamma function, therefore as a consequence we get the following theorem by similar token:

\textbf{\textsc{Theorem 2}} \textit{ Let $\Re(s)>0$ and $t\in\mathbb{R}^{+}$, then}
  \begin{equation*}
      \int_{t}^{\infty}x^{s-1}\left(\sum_{n=0}^{\infty}\phi(n)\frac{(-x)^{n}}{n!}\right)dx=\Gamma(s,t)\phi(-s).\tag{2.17}
  \end{equation*}

Theorem 1 and 2 are not as accurate as Ramanujan's Master Theorem as we see that they are not applicable to all the examples where Ramanujan's Master Theorem works. The reason to this is still unknown. But as far as we know, these Theorems are first to account for the case of finite Mellin transform.

\begin{center}
    {\large\textsc{\textbf{\S III. Discussion and Conclusion }}}
\end{center}
Ramanujan's method of deriving his theorem, which we have applied here, is quite unconventional, and this approach restricts the use of Theorem 1 to a more general class of functions. We have been only able to give limited examples so far, but with accurate numerical confrontation, thus we expect a direction in which reasonable progress can be made. Therefore, the question of extending this class remains open. 

Ramanujan also gave the following beautiful generalization of Frullani's integral theorem: Let $f$ and $g$ be continuous functions of the interval $[0,\infty)$ such that
\begin{equation*}
    f(x)-f(\infty)=\sum_{n=0}^{\infty}\frac{u
    (n)(-x)^n}{n!}\,\,\mathrm{and}\,\, g(x)-g(\infty)=\sum_{n=0}^{\infty}\frac{v
    (n)(-x)^n}{n!}\tag{3.1}
\end{equation*}
where $f(\infty)=\lim_{x\to\infty}f(x)$ and $g(\infty)=\lim_{x\to\infty}g(x)$. Assume that $u(s)/\Gamma(s+1)$ and $v(s)/\Gamma(s+1)$ satisfy the hypothesis of Hardy's version of Ramanujan's theorem. Furthermore assume that $f(0)=g(0)$ and $f(\infty)=g(\infty)$. Then for $a,b>0$ we have
\begin{equation*}
    \lim_{n\to 0^{+}}\int_{0}^{\infty}x^{n-1}(f(ax)-g(bx))dx=\left\{f(0)-f(\infty)\right\}\left\{\log{\left(\frac{b}{a}\right)+\frac{d}{ds}\left(\log\left(\frac{v(s)}{u(s)}\right)\right)_{s=0}}\right\}.\tag{3.2}
\end{equation*}
One can find the proof of the above theorem in [\cite{1}, pg. 313]. It is still unclear how one would derive the analog of Frullani's theorem for the case of finite limits using Theorem 1. What makes the proof of Eqn. (2.11) is that the limit $n\to 0^{+}$ taken on right hand side of Eqn. (1.2) is defined if we use a reflection formula for the gamma function; however, no known reflection formula is known for the lower incomplete gamma function.

Multi-dimensional generalization of Ramanujan's Master Theorem is used in evaluating loop amplitudes in Feynman diagrams \cite{8}-\cite{13}, and we believe the finite version of the Mellin transform may find applications in particle physics in some instances.

\end{document}